\documentclass[12pt]{amsart}
\usepackage{a4}
\usepackage{amssymb}

\newcommand{\R}{{\mathbb R}}
\newcommand{\C}{{\mathbb C}}

\newcommand{\LA}{{\mathfrak a}}

\newcommand{\LN}{{\mathfrak n}}

\newcommand{\LS}{{\mathfrak s}}
\newcommand{\LU}{{\mathfrak u}}
\newcommand{\LV}{{\mathfrak v}}
\newcommand{\LZ}{{\mathfrak z}}

\newcommand{\SU}{\mathrm{SU}} 
\newcommand{\U}{\mathrm{U}}

\newcommand{\Span}{\mathrm{span}} 
 
\newcommand{\RH}{\R \mathrm{H}} 
\newcommand{\CH}{\C \mathrm{H}}

\newcommand{\Max}{\mathrm{max} \, }
\newcommand{\Min}{\mathrm{min} \, }

\newcommand{\id}{\mathrm{id}}

\newcommand{\Ric}{\mathrm{Ric}}

\newcommand{\SC}{\mathrm{sc}}

\newcommand{\inner}[2]{\langle #1 , #2 \rangle} 

\newtheorem{theorem}{Theorem}[section]
\newtheorem{Prop}[theorem]{Proposition}
\newtheorem{Thm}[theorem]{Theorem}
\newtheorem{Lem}[theorem]{Lemma}
\newtheorem{Cor}[theorem]{Corollary}

\theoremstyle{definition}

\newtheorem{Def}[theorem]{Definition}

\theoremstyle{remark}

\numberwithin{equation}{section}

\begin{document}

\title[Curvatures of Lie hypersurfaces in $\CH^n$]
{Curvature properties of Lie hypersurfaces in the complex hyperbolic space}

\author{Tatsuyoshi Hamada}
\address{(Hamada) Department of Applied Mathematics, Fukuoka University, 
Fukuoka 814-0180, Japan.}
\address{(Hamada) 
JST, CREST, 5 Sanbancho, Chiyoda-ku, Tokyo 102-0075, Japan.}
\email{hamada@holst.sm.fukuoka-u.ac.jp}
\thanks{The first author was supported in part by 
Grant-in-Aid for Scientific Research (C) 21540103, 
and the third author was supported in part by 
Grant-in-Aid for Young Scientists (B) 20740040, 
The Ministry of Education, Culture, Sports, Science and Technology, Japan.}

\author{Yuji Hoshikawa}
\address{(Hoshikawa) Sakaide High School, Sakaide, Kagawa 762-0031, Japan.}

\author{Hiroshi Tamaru}
\address{(Tamaru) Department of Mathematics, Hiroshima University, 
Higashi-Hiroshima 739-8526, Japan.}
\email{tamaru@math.sci.hiroshima-u.ac.jp}


\begin{abstract}
A Lie hypersurface in the complex hyperbolic space 
is a homogeneous real hypersurface without focal submanifolds. 
The set of all Lie hypersurfaces in the complex hyperbolic space 
is bijective to a closed interval, 
which gives a deformation of homogeneous hypersurfaces 
from the ruled minimal one to the horosphere. 
In this paper, we study intrinsic geometry of Lie hypersurfaces, 
such as Ricci curvatures, scalar curvatures, and sectional curvatures. 
\end{abstract}


\medskip 

\maketitle

\section{Introduction}

The complex hyperbolic space $\CH^n$ is 
a connected and simply-connected K\"ahler manifold 
with negative constant holomorphic sectional curvature. 
A hypersurface in $\CH^n$ is said to be \textit{homogeneous} 
if it is an orbit of a closed subgroup of the isometry group of $\CH^n$. 
A homogeneous hypersurface is called a \textit{Lie hypersurface} 
if it has no focal manifolds. 
By hypersurface, we always mean a connected real hypersurface. 
The purpose of this paper is to study 
intrinsic geometry of Lie hypersurfaces in $\CH^n$, 
such as Ricci curvatures, scalar curvatures, and sectional curvatures. 
We calculate these curvatures explicitly, in terms of solvable Lie algebras. 
Our results include some interesting features of curvature properties, 
for example, 
for cases $n=2$ and $n>2$, 
the maximum values of the sectional curvatures of Lie hypersurfaces are different. 

\medskip 

One of the motivations of our study comes from submanifold geometry. 
A hypersurface in the complex hyperbolic space $\CH^n$ 
has been studied actively, 
and it provides rich sources of submanifold geometry 
(for example, refer to \cite{NR} and references therein). 
Typical examples of real hypersurfaces are given by homogeneous ones, 
which have been classified by Berndt and the third author (\cite{BT3}). 
A hypersurface in $\CH^n$, where $n \geq 2$, is homogeneous 
if and only if it is congruent to one of the following hypersurfaces: 
\begin{enumerate}
\item[(A)]
a tube around the totally geodesic $\CH^k$ for some $k = 0, \ldots , n-1$, 
\item[(B)]
a tube around the totally geodesic $\RH^n$,
\item[(N)]
a horosphere, 
\item[(S)]
the homogeneous ruled minimal hypersurface or its equidistant hypersurfaces, 
\item[(W)]
a tube around the minimal ruled submanifold $W^{2n-k}_{\varphi}$. 
\end{enumerate}
For the hypersurfaces of type (A), (B) and (N), 
their extrinsic and intrinsic geometry seem to be well investigated
(see \cite{NR}). 
Lie hypersurfaces are precisely the hypersurfaces of type (N) and (S). 
Berndt (\cite{B}) studied extrinsic geometry of Lie hypersurfaces in detail, 
which we review in Section \ref{section:extrinsic}. 
For the hypersurfaces of type (W), 
Berndt and D\'iaz-Ramos (\cite{BD}) recently studied their extrinsic geometry. 
The subject of this paper is intrinsic geometry of Lie hypersurfaces, 
which we need to understand for further study of hypersurfaces in $\CH^n$. 

\medskip 

Another motivation of our study of Lie hypersurfaces 
comes from geometry of Lie groups with left-invariant metrics. 
Lie hypersurfaces can be identified with solvable Lie groups with left-invariant metrics, 
which are called \textit{solvmanifolds}. 
In \cite{T08}, the third author constructed many Einstein solvmanifolds, 
which are homogeneous submanifolds in symmetric spaces of noncompact type. 
Searching more Einstein submanifolds, 
and finding a condition for submanifolds to be Einstein, 
are natural questions. 
Although no hypersurfaces in $\CH^n$ are Einstein 
(see \cite{R}), 
our study will be a good example of the study of intrinsic geometry of submanifolds. 

\medskip 

Lie hypersurfaces are also interesting from a view point of 
degeneration of Lie algebra (we refer to \cite{L}). 
The set of Lie hypersurfaces gives a degeneration of a solvable Lie algebra 
to nilpotent one. 
We calculate curvatures of all Lie hypersurfaces, 
which can be regarded as a study of curvature behavior under this degeneration. 

\medskip 

This paper is organized as follows. 
We recall the solvable model of the complex hyperbolic space $\CH^n$ in Section 2, 
and the Lie hypersurfaces in Section 3. 
Extrinsic geometry of Lie hypersurfaces, 
such as principal curvatures and mean curvatures, 
will be mentioned in Section 4. 
Intrinsic geometry will be studied in the remaining sections; 
We study the Ricci curvatures in Section 5, 
the scalar curvatures in Section 6, 
and the sectional curvatures in Section 7. 

%

\section{The complex hyperbolic space}

In this section, we recall the solvable model of the complex hyperbolic space 
$\CH^n$
with constant holomorphic sectional curvature $-1$. 

\begin{Def}
\label{def:solvable}
We call $(\LS, \inner{}{}, J)$ the 
\textit{solvable model of the complex hyperbolic space} if 
\begin{enumerate}
\item
$\LS$ is a Lie algebra, and there is a basis 
$\{ A_0, X_1, Y_1, \ldots X_{n-1}, Y_{n-1}, Z_0 \}$ 
whose bracket products are given by 
\begin{align*}
[A_0, X_i] = (1/2) X_i, \ 
[A_0, Y_j] = (1/2) Y_j, \ 
[A_0,Z_0] = Z_0 , \ 
[X_i, Y_i] = Z_0 , 
\end{align*}
\item
$\inner{}{}$ is an inner product on $\LS$ so that the above basis is orthonormal, 
\item
$J$ is a complex structure on $\LS$ given by 
\begin{align*}
J(A_0) = Z_0 , \ 
J(Z_0) = -A_0 , \ 
J(X_i) = Y_i , \ 
J(Y_i) = - X_i . 
\end{align*}
\end{enumerate}
\end{Def}

Throughout this paper, 
we identify our solvable model $(\LS, \inner{}{}, J)$ 
with the connected and simply-connected Lie group $S$ with Lie algebra $\LS$, 
endowed with the induced left-invariant Riemannian metric and complex structure. 

\medskip 

We note that the solvable model is a symmetric space of noncompact type
(we refer to \cite{Hel} for symmetric spaces). 
One knows 
\begin{align*}
\CH^n = \SU(1,n) / \mathrm{S}(\U(1) \times \U(n)) , 
\end{align*}
the expression of the complex hyperbolic space as a homogeneous space. 
The solvable model is nothing but the solvable part of the Iwasawa decomposition 
of $\SU(1,n)$, which can naturally be identified with $\CH^n$. 

\medskip 

We also note that $(\LS, \inner{}{})$ is a Damek-Ricci Lie algebra 
(we refer to \cite{BTV} for Damek-Ricci spaces and Lie algebras). 
Let us consider the orthogonal decomposition 
\begin{align}
\label{eq:avz}
\begin{split}
&
\LS = \LA \oplus \LV \oplus \LZ , 
\\
&
\mbox{where} \ \ 
\LA = \Span \{ A_0 \} , \ 
\LV = \Span \{ X_1, Y_1, \ldots X_{n-1}, Y_{n-1} \} , \ 
\LZ = \Span \{ Z_0 \} . 
\end{split}
\end{align}
One can directly see by definition that 
\begin{align}
\label{eq:bracket-v-w}
[V,W] = \inner{JV}{W} Z_0 \quad (\forall V, W \in \LV) . 
\end{align}
The derived subalgebra $\LN := [\LS, \LS] = \LV \oplus \LZ$ 
is the Heisenberg Lie algebra, 
which is of $H$-type from (\ref{eq:bracket-v-w}). 

\medskip 

We note that $(\LS, J)$ is a normal $j$-algebra 
(we refer to \cite{PS} for normal $j$-algebras). 
In particular, it is easy to see by definition that 
\begin{align}
\label{eq:j-inv}
\inner{JX}{JY} = \inner{X}{Y}
\quad (\forall X,Y \in \LS) , 
\end{align}
which we use in the following calculations. 

\medskip 

We here study curvature properties of the solvable model, 
which have been well-known 
(for example, 
the curvature properties of Damek-Ricci spaces can be found in \cite{BTV}). 
Let $X,Y \in \LS$. 
As in \cite{Mi}, the Kozsul formula yields that the Levi-Civita connection 
$\nabla^{\LS}$ of $(\LS, \inner{}{})$ is given by 
\begin{align}
\nabla^{\LS}_X Y = (1/2) [X,Y] + U^{\LS}(X,Y) , 
\end{align}
where $U^{\LS} : \LS \times \LS \to \LS$ is the symmetric bilinear form defined by 
\begin{align}
\label{eq:u}
2 \inner{U^{\LS}(X,Y)}{Z} = \inner{[Z,X]}{Y} + \inner{X}{[Z,Y]} \quad 
(\forall Z \in \LS) . 
\end{align}

\begin{Lem}
\label{lem:LC}
Write $X = a_1 A_0 + V + a_2 Z_0$ and $Y = b_1 A_0 + W + b_2 Z_0$ 
according to the decomposition (\ref{eq:avz}), 
where $V, W \in \LV$. 
Thus, one has 
\begin{align*}
2 \nabla^{\LS}_{X} Y = 
(\inner{V}{W} + 2 a_2 b_2) A_0 
- b_1 V - a_2 J W - b_2 J V  
+ (\inner{JV}{W} - a_2 b_1) Z_0 . 
\end{align*}
\end{Lem}

\begin{proof}
We have only to calculate $[X,Y]$ and $U^{\LS}(X,Y)$. 
Fist of all, one has 
\begin{align*}
[X,Y] = (1/2) a_1 W - (1/2) b_1 V + (\inner{JV}{W} + a_1 b_2 - a_2 b_1) Z_0 . 
\end{align*}
For $U^{\LS}(X,Y)$, we calculate its each component. 
Let $V' \in \LV$. 
One can see from (\ref{eq:bracket-v-w}) and (\ref{eq:j-inv}) that 
\begin{align*}
2 \inner{U^{\LS}(X,Y)}{V'} 
& = \inner{[V',X]}{Y} + \inner{X}{[V',Y]} \\
& = - (1/2) a_1 \inner{V'}{W} - b_2 \inner{V'}{JV} 
- (1/2) b_1 \inner{V}{V'} - a_2 \inner{V'}{JW} . 
\end{align*}
This yields that the $\LV$-component of $U^{\LS}(X,Y)$ satisfies 
\begin{align*}
2 U^{\LS}(X,Y)_{\LV} 
= - (1/2) a_1 W - b_2 JV - (1/2) b_1 V - a_2 JW . 
\end{align*}
Other components are easy to calculate, 
and hence we can calculate $\nabla^{\LS}_X Y$. 
\end{proof}

The \textit{Riemannian curvature} $R^{\LS}$ of the solvable model is defined by 
\begin{align*}
R^{\LS}(X,Y) 
:= \nabla^{\LS}_{[X,Y]} - [\nabla^{\LS}_X, \nabla^{\LS}_Y] 
= \nabla^{\LS}_{[X,Y]} 
- \nabla^{\LS}_X \nabla^{\LS}_Y + \nabla^{\LS}_Y \nabla^{\LS}_X , 
\end{align*}
which has a very simple expression in this case. 

\begin{Lem}
\label{lem:RC}
The Riemannian curvature $R^{\LS}$ satisfies 
\begin{align*}
4 R^{\LS}(X,Y)Z
= \inner{Y}{Z} X 
- \inner{X}{Z} Y 
+ \inner{JY}{Z} JX 
- \inner{JX}{Z} JY 
- 2 \inner{JX}{Y} JZ . 
\end{align*}
\end{Lem}

\begin{proof}
For the calculation, let 
\begin{align*}
X = a_1 A_0 + V + a_2 Z_0 , \ 
Y = b_1 A_0 + W + b_2 Z_0 , \ 
Z = c_1 A_0 + U + c_2 Z_0 , 
\end{align*}
where $V, W, U \in \LV$. 
It is convenient to use 
\begin{align*}
2 [X, Y] = a_1 W - b_1 V + 2 \inner{JX}{Y} Z_0 . 
\end{align*}
It follows from Lemma \ref{lem:LC} and direct calculations that 
\begin{align*}
4 (\nabla^{\LS}_{[X, Y]} Z)_{\LA} 
= & 
(a_1 \inner{W}{U} - b_1 \inner{V}{W} + 4 c_2 \inner{JX}{Y}) A_0 , \\
4 (\nabla^{\LS}_{[X, Y]} Z)_{\LV} 
= & 
b_1 c_1 V 
+ b_1 c_2 JV 
- a_1 c_1 W  
- a_1 c_2 JW 
- 2 \inner{JX}{Y} JU , \\
4 (\nabla^{\LS}_{[X, Y]} Z)_{\LZ} 
= & 
(a_1 \inner{JW}{U} 
- b_1 \inner{JV}{U} 
- 4 c_1 \inner{JX}{Y}) Z_0 , 
\end{align*}
where subscript $\LU$ denotes the $\LU$-component. 
One can also directly calculate each component of $\nabla_X \nabla_Y Z$ as 
\begin{align*}
4 (\nabla^{\LS}_X \nabla^{\LS}_Y Z)_{\LA}
= & 
(- c_1 \inner{V}{W}
- c_2 \inner{V}{JW}
- b_2 \inner{V}{JU} 
+ 2 a_2 \inner{JW}{U}
-4 a_2 b_2 c_1) A_0 , \\
4 (\nabla^{\LS}_X \nabla^{\LS}_Y Z)_{\LV} 
= & 
(- \inner{W}{U} - 2 b_2 c_2) V 
+ (\inner{W}{JU} + 2 b_2 c_1) JV \\
& \quad 
+ a_2 c_1 JW - a_2 c_2 W - a_2 b_2 U , \\
4 (\nabla^{\LS}_X \nabla^{\LS}_Y Z)_{\LZ}
= & 
(- c_1 \inner{JV}{W} 
- c_2 \inner{V}{W} 
- b_2 \inner{V}{U}
- 2 a_2 \inner{W}{U} 
- 4 a_2 b_2 c_2) Z_0 . 
\end{align*}
One can obtain $\nabla_Y \nabla_X Z$ by symmetry. 
By summing up them, we complete the proof. 
\end{proof}

The \textit{Ricci operator} $\Ric^{\LS}$ is defined by 
\begin{align*}
\Ric^{\LS}(X) := \sum R^{\LS}(E_i, X) E_i , 
\end{align*}
where $\{ E_i \}$ is an orthonormal basis of $\LS$. 
A Riemannian manifold is said to be \textit{Einstein} 
if the Ricci operator is a scalar map, $\Ric = c \cdot \id$, 
and the scalar $c$ is called the \textit{Einstein constant}. 

\begin{Prop}
$\Ric^{\LS}(X) = - ((n+1)/2) X$ holds for every $X \in \LS$. 
Hence, the solvable model is an Einstein manifold with negative Einstein constant. 
\end{Prop}

\begin{proof}
First of all, it follows from Lemma \ref{lem:RC} that, for $X, Y \in \LS$, 
\begin{align}
\label{eq:RC=3/4}
R^{\LS}(X,Y)X
= - (3/4) \inner{JX}{Y} JX - (1/4) |X|^2 Y + (1/4) \inner{X}{Y} X . 
\end{align}
Let $X \in \LS$. 
We may and do assume that $|X|=1$ without loss of generality. 
We use an orthonormal basis $\{ E_i \}$ such that $E_1 = X$ and $E_2 = JX$. 
Thus, one can see from (\ref{eq:RC=3/4}) that 
\begin{align*}
R^{\LS}(E_1, X) E_1 = 0 , \
R^{\LS}(E_2, X) E_2 = -X , \ 
R^{\LS}(E_k, X) E_k = -(1/4) X
\end{align*}
for $k = 3, \ldots, 2n$. 
By adding them, we obtain the Ricci operator. 
\end{proof}

The \textit{sectional curvature} $K^{\LS}_{\sigma}$ 
of a plane $\sigma$ in $\LS$ is defined by 
\begin{align*}
K^{\LS}_{\sigma} := \inner{R^{\LS}(X,Y)X}{Y} , 
\end{align*}
where $\{ X, Y \}$ is an orthonormal basis of $\sigma$. 
Recall that the K\"ahler angle $\alpha$ of $\sigma$ is given by 
$\cos(\alpha) = \inner{JX}{Y}$. 
The holomorphic sectional curvature is the sectional curvature 
of a complex plane, that is, a plane with $\alpha = 0$. 

\begin{Prop}
\label{prop:ch-n}
$K^{\LS}_{\sigma} = - (1/4) - (3/4) \inner{JX}{Y}^2$ holds 
for a plane $\sigma$ with an orthonormal basis $\{ X, Y \}$. 
Hence, the solvable model has the 
constant holomorphic sectional curvature $-1$. 
\end{Prop}

\begin{proof}
It is direct from (\ref{eq:RC=3/4}) 
and the definition of the holomorphic sectional curvature. 
\end{proof}

One can see that all the curvatures can be calculated 
in terms of the Lie algebra $\LS$. 
The curvatures of Lie hypersurfaces are more complicated, 
but they can also be completely calculated in terms of the Lie algebras. 

\section{Lie hypersurfaces}
\label{lie-hypersurface}

The Lie hypersurfaces have been introduced and studied by Berndt (\cite{B}). 
In this section, we recall his results on Lie hypersurfaces, 
and mention a calculation of 
the second fundamental forms and the Levi-Civita connections. 

\begin{Def}
A homogeneous hypersurface in the complex hyperbolic space is called a 
\textit{Lie hypersurface} 
if it has no focal submanifolds. 
\end{Def}

Note that this definition looks different 
from the original definition by Berndt (\cite{B}). 
In \cite{B}, a Lie hypersurface is defined by an orbit of 
a codimension one subgroup $S'$ of $S$, 
where $S$ is the solvable Lie group involved in the solvable model. 
It is obvious that every orbit of $S'$ is a hypersurface. 
Therefore, every orbit of $S'$ is a Lie hypersurface in our sense. 
The converse of this statement also holds. 
It follows from the following; 
every cohomogeneity one action without singular orbit 
(that is, an isometric action all of whose orbits are codimension one) 
is orbit equivalent to an action of a codimension one subgroup of $S$. 
In fact, this is true not only for $\CH^n$, 
but also for every symmetric space of noncompact type 
(\cite{BT1}). 

\begin{Thm}[\cite{B}]
\label{thm:classification}
Every Lie hypersurface in the complex hyperbolic space is 
isometrically congruent to the orbit $S(\theta).o$ for $\theta \in [0, \pi/2]$, 
where $o$ is the origin and 
$S(\theta)$ is the connected Lie subgroup of $S$ with Lie algebra 
\begin{align*}
\LS(\theta) := \LS \ominus \R (\cos(\theta) X_1 + \sin(\theta) A_0) . 
\end{align*}
\end{Thm}

Note that $\ominus$ denotes the orthogonal complement. 
One can check that $\LS(\theta)$ is a subalgebra of $\LS$, 
which is obviously of codimension one. 
We will use the unit normal vector 
\[
\xi := \cos(\theta) X_1 + \sin(\theta) A_0 
\]
and the orthogonal decomposition 
\begin{align}
\label{eq:v0}
\begin{split}
&
\LS(\theta) = \Span \{ T \} \oplus \Span \{ Y_1 \} \oplus \LV_0 \oplus \LZ , 
\\
&
\mbox{where} \ \ 
T := \cos(\theta) A_0 - \sin(\theta) X_1 , \ \ 
\LV_0 := \Span \{ X_2, Y_2, \ldots, X_{n-1}, Y_{n-1} \} . 
\end{split}
\end{align}
We use this decomposition frequently. 

\medskip 

Note that the set of Lie hypersurfaces gives a degeneration of Lie algebra, 
from a solvable Lie algebra $\LS(0)$ to a nilpotent Lie algebra $\LS(\pi/2)$. 
If $\theta \neq \pi/2$, 
then $\LS(\theta)$ is solvable but not nilpotent. 
The Lie algebra $\LS(\pi/2) = \LN$ is the Heisenberg Lie algebra, 
which is nilpotent. 
Our study of this paper can be regarded as a study of curvature behavior 
under a degeneration of Lie algebra. 

\medskip 

The orbit $S(0).o$ is the homogeneous ruled minimal hypersurface 
(we refer to \cite{LR}), 
and the orbit $S(\pi/2).o$ is a horosphere. 
An equidistant hypersurface to $S(0).o$ is an orbit $S(0).\gamma(t)$ for $t>0$, 
where $\gamma$ is a normal geodesic of $S(0).o$ starting from $o$. 
The left-translation by $\gamma(t)^{-1} \in S$ gives 
\begin{align*}
S(0).\gamma(t) 
\cong (\gamma(t)^{-1} S(0) \gamma(t)).o 
= S(\theta).o 
\end{align*}
for $\theta \in (0, \pi/2)$. 
Thus, the set of Lie hypersurfaces gives a deformation of hypersurfaces, 
keeping homogeneity, 
from the ruled minimal homogeneous hypersurface $S(0).o$ 
to the horosphere $S(\pi/2).o$. 
We refer to \cite{B} and \cite{BT1} for detail. 

\medskip 

In the rest of the paper, we denote by $S(\theta)$ instead of $S(\theta).o$. 
We need the Levi-Civita connections $\nabla$
and the second fundamental forms $h$ of the Lie hypersurfaces $S(\theta)$. 
One knows the formula 
\begin{align}
\label{eq:h}
\nabla^{\LS}_X Y = \nabla_X Y + h(X,Y) 
\end{align}
for $X, Y \in \LS(\theta)$, 
where $\nabla_X Y \in \LS(\theta)$ and $h(X,Y) \in \R \xi$. 

\begin{Prop}
\label{prop:h}
Let $X, Y \in \LS(\theta)$ and write as 
$X = a_1 T + a_2 Y_1 + V + a_3 Z_0$ 
and 
$Y = b_1 T + b_2 Y_1 + W + b_3 Z_0$ 
for $V,W \in \LV_0$ according to the decomposition (\ref{eq:v0}). 
Thus, the second fundamental form $h$ satisfies 
\begin{align*}
2 h(X,Y)
= 
\left(
(\inner{X}{Y} + a_3 b_3) \sin(\theta)
+ (a_2 b_3 + a_3 b_2) \cos(\theta) 
\right) \xi . 
\end{align*}
\end{Prop}

\begin{proof}
We know $\nabla^{\LS}$ by Lemma \ref{lem:LC}. 
Thus the second fundamental form $h$ can be calculated directly from 
$h(X,Y) = \inner{\nabla^{\LS}_X Y}{\xi} \xi$. 
\end{proof}

The Levi-Civita connection $\nabla$ can be calculated in terms of $h$. 

\begin{Prop}
\label{prop:nabla}
For the Lie hypersurface $S(\theta)$, we have 
\begin{enumerate}
\item
$\nabla_{T} T = 0$, \ 
$\nabla_{T} Y_1  = - (1/2) \sin(\theta) Z_0$, \ 
$\nabla_{T} V  = 0$, \ 
$\nabla_{T} Z  = (1/2) \sin(\theta) Y_1$, 
\item
$\nabla_{Y_1} T = - (1/2) \cos(\theta) Y_1 + (1/2) \sin(\theta) Z_0$, \\
$\nabla_{Y_1} Y_1 = (1/2) \cos(\theta) T$, \ 
$\nabla_{Y_1} V = 0$, \ 
$\nabla_{Y_1} Z_0 = - (1/2) \sin(\theta) T$,
\item
$\nabla_{V} T = - (1/2) \cos(\theta) V$, \ 
$\nabla_{V} Y_1 = 0$, \\ 
$\nabla_{V} W = (1/2) [V,W] + (1/2) \cos(\theta) \inner{V}{W} T$, \ 
$\nabla_{V} Z_0 = - (1/2) J V$, 
\item
$\nabla_{Z_0} T = (1/2) \sin(\theta) Y_1 - \cos(\theta) Z_0$, \ 
$\nabla_{Z_0} Y_1 = - (1/2) \sin(\theta) T$, \\ 
$\nabla_{Z_0} W = - (1/2) J W$, \ 
$\nabla_{Z_0} Z_0 = \cos(\theta) T$. 
\end{enumerate}
\end{Prop}

\begin{proof}
It is direct from 
$\nabla_X Y = \nabla^{\LS}_X Y - h(X,Y)$. 
\end{proof}

\section{Extrinsic geometry}
\label{section:extrinsic}

Extrinsic geometry of the Lie hypersurfaces $S(\theta)$ have been studied 
by Berndt (\cite{B}) in detail. 
In this section we recall his results, 
since the proofs follow easily from Proposition \ref{prop:h}, 
and the calculations are similar to the ones in the following sections. 

\medskip 

First of all we calculate the shape operator of $S(\theta)$. 
The shape operator $A_{\xi} : \LS(\theta) \to \LS(\theta)$ is determined, 
in terms of the second fundamental form $h$, by 
\begin{align}
\label{eq:shape}
\inner{A_{\xi}(X)}{Y} = \inner{h(X,Y)}{\xi} . 
\end{align}

\begin{Thm}
\label{thm:shape-op}
The shape operator of the Lie hypersurface $S(\theta)$ satisfies 
\begin{enumerate}
\item
$A_{\xi}(T) = (1/2) \sin(\theta) T$, 
\item
$A_{\xi}(Y_1) = (1/2) \sin(\theta) Y_1 + (1/2) \cos(\theta) Z_0$, 
\item
$A_{\xi}(V) = (1/2) \sin(\theta) V$ \ for $V \in \LV_0$, 
\item
$A_{\xi}(Z_0) = (1/2) \cos(\theta) Y_1 + \sin(\theta) Z_0$. 
\end{enumerate}
\end{Thm}

\begin{proof}
One has from (\ref{eq:shape}) that 
\begin{align*}
\textstyle 
A_{\xi}(X) 
= \sum \inner{A_{\xi}(X)}{E_i} E_i 
= \sum \inner{h(X,E_i)}{\xi} E_i , 
\end{align*}
where $\{ E_i \}$ is an orthonormal basis of $\LS(\theta)$. 
We only demonstrate $A_{\xi}(V)$ for $V \in \LV_0$. 
We can assume without loss of generality that 
$|V|=1$ and the orthonormal basis $\{ E_i \}$ satisfies $E_1 = V$. 
One knows $h$ by Proposition \ref{prop:h}. 
If $i \ge 2$, then $\inner{V}{E_i} = 0$ and hence $h(V,E_i) = 0$ holds. 
Therefore, one has 
\begin{align*}
A_{\xi}(V) 
= \inner{h(V,V)}{\xi} V 
= (1/2) \sin(\theta) V , 
\end{align*}
which completes (3). 
One can prove (1), (2) and (4) by similar calculations. 
\end{proof}

An eigenvalue of the Shape operator $A_{\xi}$ is called 
a \textit{principal curvature}, 
and the dimension of the eigenspace is called the \textit{multiplicity}. 

\begin{Cor}
\label{cor:principal-curvature}
The principal curvatures of the Lie hypersurface $S(\theta)$ are 
$\lambda_1$, $\lambda_2$ and $\lambda_3$, where 
\begin{enumerate}
\item[]
$\lambda_1 := (3/4) \sin(\theta) - (1/4) \left( 1+3 \cos^2(\theta) \right)^{1/2}$, 
\item[]
$\lambda_2 := (1/2) \sin(\theta)$, 
\item[]
$\lambda_3 := (3/4) \sin(\theta) + (1/4) \left( 1+3 \cos^2(\theta) \right)^{1/2}$. 
\end{enumerate}
More precisely, 
\begin{enumerate}
\item
If $\theta \neq \pi/2$, 
then $\lambda_1 < \lambda_2 < \lambda_3$ holds 
and the multiplicities are $1$, $2n-3$ and $1$, respectively. 
\item
If $\theta = \pi/2$, 
then $\lambda_1 = \lambda_2 = 1/2 < \lambda_3 = 3/4$ holds 
and the multiplicities are $2n-2$ and $1$, respectively. 
\end{enumerate}
\end{Cor}

\begin{proof}
It follows easily from Theorem \ref{thm:shape-op} that 
$\lambda_1$, $\lambda_2$, $\lambda_3$ are principal curvatures. 
Furthermore, one can see that 
\begin{align*}
4(\lambda_3 - \lambda_2) 
= \sin(\theta) + \left( 1+3 \cos^2(\theta) \right)^{1/2} > 1 > 0 . 
\end{align*}
For convenience, let 
$\rho := \left( 1+3 \cos^2(\theta) \right)^{1/2} + \sin(\theta) > 0$. 
Thus, one can see that 
\begin{align*}
4 (\lambda_2 - \lambda_1) \rho 
= ( \left( 1+3 \cos^2(\theta) \right)^{1/2} - \sin(\theta) ) \rho 
= 4 \cos^2(\theta)
\ge 0 . 
\end{align*}
This yields $\lambda_2 \ge \lambda_1$, 
and the equality holds if and only if $\theta = \pi/2$. 
\end{proof}

The mean value of the principal curvatures is called the \textit{mean curvature}. 
A submanifold is called \textit{austere} if the set of principal curvatures 
counted with multiplicities are invariant by $-1$, 
and is called \textit{minimal} if the mean curvature is $0$. 
Obviously, an austere submanifold is minimal. 

\begin{Cor}
The mean curvature of the Lie hypersurface $S(\theta)$ is 
\begin{align*}
(n / (2n-1)) \sin (\theta) , 
\end{align*}
and hence $S(\theta)$ is minimal if and only if $\theta = 0$. 
Furthermore, $S(0)$ is an austere submanifold. 
\end{Cor}

\begin{proof}
It is easy from 
Corollary \ref{cor:principal-curvature} 
that the mean curvature is $(n / (2n-1)) \sin (\theta)$. 
Hence $S(\theta)$ is minimal if and only if $\theta = 0$. 
If $\theta = 0$, then the principal curvatures are 
$\lambda_1 = - 1/2$, $\lambda_2 = 0$, $\lambda_3 = 1/2$, 
with multiplicities $1$, $2n-3$, $1$, respectively. 
Thus $S(0)$ is austere. 
\end{proof}

The tangent vector field $J \xi$ is called the \textit{structure vector field}. 
A hypersurface in $\CH^n$ is said to be \textit{Hopf} 
if the structure vector field is an eigenvector of the shape operator. 

\begin{Cor}
\label{cor:hopf}
The Lie hypersurface $S(\theta)$ is Hopf if and only if $\theta = \pi/2$. 
\end{Cor}

\begin{proof}
By definition, the structure vector field is 
$J \xi = \cos(\theta) Y_1 + \sin(\theta) Z_0$. 
Thus, Theorem \ref{thm:shape-op} yields that 
\[
A_{\xi}(J \xi) 
= \sin(\theta) \cos(\theta) Y_1 + ((1/2) \cos^2(\theta) + \sin^2(\theta)) Z_0 . 
\]
Since $J \xi, A_{\xi} (J \xi) \in \Span \{ Y_1, Z \}$, 
one has that $J \xi$ is an eigenvector if and only if 
\begin{align*}
0 = \inner{A_{\xi}(J \xi)}{\sin(\theta) Y_1 - \cos(\theta) Z_0} 
= - (1/2) \cos^3(\theta) . 
\end{align*}
This concludes the claim. 
\end{proof}

\section{Ricci curvatures}

In this section, we calculate the Ricci curvatures of the Lie hypersurfaces 
$S(\theta)$. 
Let $R$ be the Riemannian curvature of $S(\theta)$, 
and recall that the Ricci operator 
$\Ric : \LS(\theta) \to \LS(\theta)$ 
is defined by 
\begin{align*}
\textstyle \Ric(X) = \sum R(E_i, X) E_i , 
\end{align*}
where $\{ E_i \}$ is an orthonormal basis of $\LS(\theta)$. 
We use the orthogonal decomposition 
$\LS(\theta) = \R T \oplus \R Y_1 \oplus \LV_0 \oplus \LZ$ 
given in (\ref{eq:v0}). 

\begin{Thm}
\label{thm:ricci}
The Ricci operator of the Lie hypersurface 
$(\LS(\theta), \inner{}{})$ satisfies 
\begin{enumerate}
\item
$\Ric(T) = - (1/4) (2 + (2n-1) \cos^2(\theta)) T$, 
\item
$\Ric(Y_1) = 
- (1/4) (2 + (2n-3) \cos^2(\theta)) Y_1 + (n/2) \sin(\theta) \cos(\theta) Z_0$, 
\item
$\Ric(V) = - (1/4) (2 + (2n-1) \cos^2(\theta)) V$ 
\quad for $V \in \LV_0$, 
\item
$\Ric(Z_0) = 
(n/2) \sin(\theta) \cos(\theta) Y_1 + (1/2) ((n-1) - 2n \cos^2(\theta)) Z_0$. 
\end{enumerate}
\end{Thm}

\begin{proof}
We only demonstrate $\Ric(V)$ for $V \in \LV_0$. 
By definition, one has 
\begin{align*}
\textstyle 
\Ric(V) 
= 
R(T, V)T 
+ R(Y_1, V)Y_1
+ \sum R(W_i, V)W_i
+ R(Z, V)Z , 
\end{align*}
where $\{ W_i \}$ is an orthonormal basis of $\LV_0$. 
One knows $\nabla$ by Proposition \ref{prop:nabla}, which yields that 
\begin{align*}
R(T, V) T 
& = 
(1/2) \cos(\theta) \nabla_V T , 
\\
R(Y_1, V) Y_1
& = 
(1/2) \cos(\theta) \nabla_V T , 
\\
R(W_i, V) W_i 
& = 
- \nabla_{[V,W_i]} W_i 
- (1/2) \nabla_{W_i} [V,W_i] 
\\
& \qquad 
- (1/2) \cos(\theta) \inner{V}{W_i} \nabla_{W_i} T 
+ (1/2) \cos(\theta) \nabla_V T , 
\\
R(Z_0, V) Z_0 
& = 
\cos(\theta) \nabla_V T 
+ (1/2) \nabla_{Z_0} JV . 
\end{align*}
Without loss of generality, we may and do assume that 
$|V|=1$, $W_1 = V$, and $W_2 = JV$. 
Note that $[V, JV] = Z_0$, and $[V, W_i] = 0$ for $i \neq 2$ 
(see Eq.\ (\ref{eq:bracket-v-w})). 
Therefore, one has 
\begin{align*}
\textstyle 
\sum R(W_i, V) W_i 
& = 
- \nabla_{Z_0} JV 
- (1/2) \nabla_{JV} Z_0 
\\
& \qquad 
- (1/2) \cos(\theta) \nabla_V T 
+ ((2n-4)/2) \cos(\theta) \nabla_V T . 
\end{align*}
Altogether, one get 
\begin{align*}
\Ric(V) 
& = 
((2n-1)/2) \cos(\theta) \nabla_V T 
- (1/2) \nabla_{Z_0} JV
- (1/2) \nabla_{JV} Z_0 
\\
& = 
- ((2n-1)/4) \cos^2(\theta) V 
+ (1/4) J^2 V
+ (1/4) J^2 V 
\\
& = 
- (1/4) (2 + (2n-1) \cos^2(\theta)) V . 
\end{align*}
One can calculate $\Ric(T)$, $\Ric(Y_1)$ and $\Ric(Z_0)$ in the similar way. 
\end{proof}

An eigenvalue of the Ricci operator is called a 
\textit{principal Ricci curvature}, 
and the dimension of the eigenspace is called the \textit{multiplicity}. 

\begin{Cor}
The principal Ricci curvatures of $S(\theta)$ are 
$\alpha_1$, $\alpha_2$ and $\alpha_3$, where 
\begin{align*}
\alpha_1 
& = 
(n-2)/4 
- ((6n-3)/8) \cos^2(\theta)
- (1/8) A^{1/2} , 
\\
\alpha_2 
& = 
- (1/2) 
- ((2n-1)/4) \cos^2(\theta) , 
\\
\alpha_3 
& = 
(n-2)/4
- ((6n-3)/8) \cos^2(\theta)
+ (1/8) A^{1/2} , 
\end{align*}
here $A$ is defined by 
\begin{align*}
A = 4 n^2 + 4n (2n - 3) \cos^2(\theta) -3 (2n+1)(2n-3) \cos^4(\theta) . 
\end{align*}
More precisely, 
\begin{enumerate}
\item
If $\theta \neq \pi/2$, 
then $\alpha_1 < \alpha_2 < \alpha_3$ holds and 
the multiplicities are $1$, $2n-3$ and $1$, respectively. 
\item
If $\theta = \pi/2$, 
then $\alpha_1 = \alpha_2 = - 1/2 < \alpha_3 = (n-1)/2$ holds 
and the multiplicities are $2n-2$ and $1$, respectively. 
\end{enumerate}
\end{Cor}

\begin{proof}
Theorem \ref{thm:ricci} yields that 
$T$ and $V \in \LV_0$ are eigenvectors with eigenvalue $\alpha_2$. 
Note that $\dim \R T \oplus \LV_0 = 2n-3$. 
A direct calculation yields that the eigenvalues of $\Ric$ on $\R Y_1 \oplus \LZ$ 
are $\alpha_1$ and $\alpha_3$. 
To complete the proof, we need to compare $\alpha_1$, $\alpha_2$ and $\alpha_3$. 
First of all, one has $\alpha_2 < \alpha_3$, since 
\begin{align*}
8(\alpha_3 - \alpha_2) 
= 2n + (- 2n + 1) \cos^2(\theta) + A^{1/2}
\ge 2n (1 - \cos^2(\theta)) + \cos^2(\theta) 
> 0 . 
\end{align*}
To show $\alpha_1 \le \alpha_2$, 
let $\rho = A^{1/2} + (2n + (- 2n + 1) \cos^2(\theta)) > 0$. 
Thus, one can see that 
\begin{align*}
8(\alpha_2 - \alpha_1) \rho = & 
( A^{1/2} - (2n + (- 2n + 1) \cos^2(\theta)) ) \rho \\
= & 
16n(n-1) \cos^2(\theta) (1 - \cos^2(\theta)) + 8 \cos^4(\theta) \\
\ge & 0 . 
\end{align*}
The equality holds if and only if $\cos(\theta) = 0$. 
\end{proof}

If $\theta = 0$, then $\alpha_1 < \alpha_2 < \alpha_3 < 0$ holds. 
Thus, the homogeneous ruled minimal hypersurface $S(0)$, 
and also $S(\theta)$ for small $\theta$, have negative Ricci curvatures. 
Furthermore, 
If $\theta \to \pi/2$, 
that is the Lie hypersurface $S(\theta)$ degenerates to the horosphere $S(\pi/2)$, 
then the number of distinct principal curvatures goes to two. 

\section{Scalar curvatures}

In this small section, we mention the scalar curvatures 
of Lie hypersurfaces. 
Theorem \ref{thm:ricci} immediately yields that 

\begin{Thm}
The scalar curvature $\SC$ of the Lie hypersurface 
$(\LS(\theta), \inner{}{})$ 
satisfies 
\begin{align*}
\SC = - (n-1)/2 - (n(2n-1)/2) \cos^2(\theta) . 
\end{align*}
\end{Thm}

Therefore, every Lie hypersurface has a negative scalar curvature. 

\section{Sectional curvatures}

In this section, we calculate the sectional curvatures of Lie hypersurfaces 
and determine the maximum values of them. 
It seems to be surprising that 
the maximum values for $n=2$ and $n>2$ are different. 

\begin{Thm}
\label{thm:sectional-curvature}
Let $\sigma$ be a plane in $\LS(\theta)$ with an orthonormal basis $\{ X, Y \}$. 
According to the decomposition 
$\LS(\theta) = \R T \oplus \R Y_1 \oplus \LV_0 \oplus \LZ$ 
given in (\ref{eq:v0}), we write 
\begin{eqnarray*}
X = a_1 T + a_2 Y_1 + V + a_3 Z_0 , \quad 
Y = b_1 T + b_2 Y_1 + W + b_3 Z_0 . 
\end{eqnarray*}
Thus, the sectional curvature $K_{\sigma}$ satisfies 
\begin{align*}
K_{\sigma} 
& = 
- (1/4) - (3/4) \inner{JX}{Y}^2 
+ (1/4) (1 + a_3^2 + b_3^2) \sin^2(\theta) 
\\
& \qquad 
+ (1/2) (a_2 a_3 + b_2 b_3) \sin(\theta) \cos(\theta) 
- (1/4) (a_2 b_3 - a_3 b_2)^2 \cos^2(\theta) . 
\end{align*}
\end{Thm}

\begin{proof}
Let $K^{\LS}_{\sigma}$ be the sectional curvature of $\sigma$ in $S$. 
Proposition \ref{prop:ch-n} yields that 
\begin{align*}
K_{\sigma} 
= K^{\LS}_{\sigma} + (K_{\sigma} - K^{\LS}_{\sigma}) 
= - (1/4) - (3/4) \inner{JX}{Y}^2 + (K_{\sigma} - K^{\LS}_{\sigma}) . 
\end{align*}
To calculate $K_{\sigma} - K^{\LS}_{\sigma}$, 
we recall the equation of Gauss (see \cite{KN}, Chapter VII): 
\begin{align*}
\inner{R(X,Y)Z}{W} 
= 
\inner{R^{\LS}(X,Y)Z}{W} 
+ \inner{h(X,Z)}{h(Y,W)} - \inner{h(Y,Z)}{h(X,W)} . 
\end{align*}
Note that $R$ and $R^{\LS}$ are the Riemannian curvatures of 
$\LS(\theta)$ and $\LS$ respectively, 
and $h$ is the second fundamental form. 
This immediately concludes that 
\begin{align*}
K_{\sigma} - K^{\LS}_{\sigma} 
= \inner{h(X,X)}{h(Y,Y)} - |h(X,Y)|^2 . 
\end{align*}
One knows $h$ by Proposition \ref{prop:h}, which yields that 
\begin{align*}
K_{\sigma} - K^{\LS}_{\sigma}
& = 
(1/4) (1 + a_3^2 + b_3^2) \sin^2(\theta) 
+ (1/2) (a_2 a_3 + b_2 b_3) \sin(\theta) \cos(\theta) 
\\
& \qquad 
- (1/4) (a_2 b_3 - a_3 b_2)^2 \cos^2(\theta) . 
\end{align*}
This completes the proof. 
\end{proof}

We will study the maximum and the minimum of the sectional curvature 
of each Lie hypersurface $S(\theta)$. 
At first, we study the case $n=2$. 

\begin{Cor}
Let $n=2$. 
The sectional curvature $K$ of the Lie hypersurface $S(\theta)$ 
in the complex hyperbolic plane $\CH^2$ satisfies 
\begin{align*}
\Max K_{\sigma}
= & 
- \frac{1}{4} 
- \frac{3}{8} \cos^2(\theta) 
+ \frac{1}{8} \sqrt{16 \sin^4(\theta) 
+ 9 \cos^4(\theta) + 40 \sin^2(\theta) \cos^2(\theta)} , \\
\Min K_{\sigma}
= & 
- \frac{1}{4} 
- \frac{3}{8} \cos^2(\theta) 
- \frac{1}{8} \sqrt{16 \sin^4(\theta) 
+ 9 \cos^4(\theta) + 40 \sin^2(\theta) \cos^2(\theta)} . 
\end{align*}
\end{Cor}

\begin{proof}
Since $n=2$, we have $\LV_0 = 0$ and $\LS(\theta) = \Span \{ T, Y_1, Z_0 \}$. 
Let $\sigma$ be a plane of $\LS(\theta)$. 
One can take an orthonormal basis $\{ X, Y \}$ such that 
\begin{align*}
X & 
= a_1 T + a_2 Y_1 + a_3 Z_0 
= a_1 \cos(\theta) A_0 - a_1 \sin(\theta) X_1 + a_2 Y_1 + a_3 Z_0 , \\ 
Y & 
= b_2 Y_1 + b_3 Z_0 . 
\end{align*}
By Theorem \ref{thm:sectional-curvature}, one can see that 
\begin{align*}
K_{\sigma} 
& = 
- (1/4) - (3/4) a_1^2 (b_3 \cos(\theta) - b_2 \sin(\theta))^2 
+ (1/4) (1 + a_3^2 + b_3^2) \sin^2(\theta) 
\\
& \qquad 
+ (1/2) (a_2 a_3 + b_2 b_3) \sin(\theta) \cos(\theta) 
- (1/4) (a_2 b_3 - a_3 b_2)^2 \cos^2(\theta) . 
\end{align*}
Since $\{ X, Y \}$ is orthonormal, there exists $t \in \R$ such that 
$(b_2, b_3) = (\cos(t), \sin(t))$ and 
$(a_2, a_3) = c (\sin(t), - \cos(t))$, 
where $c^2 = 1 - a_1^2$. 
A straightforward calculation yields that 
\begin{align*}
K_{\sigma} = 
- (1/4) + (1/2) \sin^2(\theta) - (1/4) \cos^2(\theta) + (1/4) a_1^2 B(t) , 
\end{align*}
where $B(t)$ is defined by 
\begin{align*}
B(t) 
= & 
\cos^2(\theta) 
- 4 \sin^2(\theta) \cos^2(t)
- 3 \cos^2(\theta) \sin^2(t) 
+ 2 \sin(\theta) \cos(\theta) \sin(t) \cos(t) \\ 
= & 
- 2 \sin^2(\theta) 
- (1/2) \cos^2(\theta)
\\
& \qquad 
+ 4 \sin(\theta) \cos(\theta) \sin(2t) 
+ (1/2) (-4 \sin^2(\theta) + 3 \cos^2(\theta)) \cos(2t) . 
\end{align*}
By applying 
$- (A^2 + B^2)^{1/2} \leq A \cos(2t) + B \sin(2t) \leq (A^2 + B^2)^{1/2}$, 
we have 
\begin{align*}
\Max B(t) 
= & 
- 2 \sin^2(\theta) 
- (1/2) \cos^2(\theta)
\\
& \qquad 
+ (1/2) 
\left(
64 \sin(\theta) \cos(\theta) + (-4 \sin^2(\theta) + 3 \cos^2(\theta))^2
\right)^{1/2} , \\ 
\Min B(t) 
= & 
- 2 \sin^2(\theta) 
- (1/2) \cos^2(\theta)
\\
& \qquad 
- (1/2) 
\left( 
64 \sin(\theta) \cos(\theta) + (-4 \sin^2(\theta) + 3 \cos^2(\theta))^2
\right)^{1/2} . 
\end{align*}
Note that $\Max B(t) \ge 0$ and $\Min B(t) \le 0$. 
Therefore, $K_{\sigma}$ attains the maximum (resp.\ minimum)
if and only if $a_1^2 = 1$ and $B(t)$ attains the maximum (resp.\ minimum). 
This finishes the proof. 
\end{proof}

We study next the case $n>2$. 

\begin{Cor}
Let $n>2$. 
Then the sectional curvature $K$ of the Lie hypersurface $S(\theta)$ 
in the complex hyperbolic space $\CH^n$ satisfies 
\begin{align*}
\Max K_{\sigma} 
= & 
- \frac{1}{4} 
+ \frac{3}{8} \sin^2(\theta)
+ \frac{1}{8} \sin(\theta) \sqrt{\sin^2(\theta) + 4 \cos^2(\theta)} . 
\end{align*}
\end{Cor}

\begin{proof}
One can see that 
\begin{align*}
K_{\sigma} 
= K^{\LS}_{\sigma} + (K_{\sigma} - K^{\LS}_{\sigma}) 
\leq - (1/4) + \Max (K_{\sigma} - K^{\LS}_{\sigma}) . 
\end{align*}
We will find a plane $\sigma$ which attains the maximums 
of $K^{\LS}_{\sigma}$ and $K_{\sigma} - K^{\LS}_{\sigma}$, simultaneously. 
Let $\sigma$ be a plane in $\LS(\theta)$ with an orthonormal basis $\{ X, Y \}$. 
We write 
\begin{eqnarray*}
X = a_1 T + a_2 Y_1 + V + a_3 Z_0 , \quad 
Y = b_1 T + b_2 Y_1 + W + b_3 Z_0 
\end{eqnarray*}
as in Theorem \ref{thm:sectional-curvature}. 
By changing an orthonormal basis of $\sigma$, 
we may and do assume that $b_3 = 0$ without loss of generality. 
Thus, Theorem \ref{thm:sectional-curvature} yields that 
\begin{align*}
K_{\sigma} 
& \leq 
- (1/4) 
+ (1/4) (1 + a_3^2) \sin^2(\theta) 
+ (1/2) a_2 a_3 \sin(\theta) \cos(\theta) 
\\
& = 
- (1/4)
+ (1/4) \sin^2(\theta) 
+ (1/8) \sin(\theta) 
\left( 2 \sin(\theta) a_3^2 + 4 \cos(\theta) a_2 a_3 \right) . 
\end{align*}
The equality holds if, for example, $\inner{JX}{Y} = 0$ and $b_2=0$. 
Since $|X|=1$, there exist $r \in [0,1]$ and $t \in \R$ such that 
$(a_2, a_3) = (r \cos(t), r \sin(t))$. 
We thus have 
\begin{align*}
2 \sin(\theta) a_3^2
+ 4 \cos(\theta) a_2 a_3 
& = 
2 \sin(\theta) r^2 \sin^2(t) 
+ 4 \cos(\theta) r^2 \sin(t) \cos(t) \\
& = 
r^2
\left(
\sin(\theta) 
- \sin(\theta) \cos(2t)
+ 2 \cos(\theta) \sin(2t)
\right) 
\\
& \leq 
r^2
\left(
\sin(\theta) 
+ (\sin^2(\theta) + 4 \cos^2(\theta))^{1/2} 
\right) 
\\
& \leq 
\sin(\theta) 
+ (\sin^2(\theta) + 4 \cos^2(\theta))^{1/2} . 
\end{align*}
This concludes that 
\begin{align}
\label{eq:n>2}
K_{\sigma} 
\le 
- (1/4)
+ (3/8) \sin^2(\theta) 
+ (1/8) \sin(\theta) (\sin^2(\theta) + 4 \cos^2(\theta))^{1/2} . 
\end{align}
The equality holds for certain $X = a_2 Y_1 + a_3 Z_0$ and $Y = V$, 
and therefore, the right hand side coincides with the maximum of $K_{\sigma}$. 
\end{proof}

Note that the inequality (\ref{eq:n>2}) holds for every $n$. 
Thus, we have $\Max K_{\sigma}^{(n>2)} \ge \Max K_{\sigma}^{(n=2)}$, 
where $K_{\sigma}^{(n>2)}$ and $K_{\sigma}^{(n=2)}$ 
are the sectional curvatures for cases $n>2$ and $n=2$, respectively. 
It is natural to ask when the equality holds. 

\begin{Prop}
We have $\Max K_{\sigma}^{(n>2)} \ge \Max K_{\sigma}^{(n=2)}$, 
and the equality holds if and only if $\theta = 0, \pi/2$. 
\end{Prop}

\begin{proof}
For convenience, let us define 
\begin{align*}
C & := 3 + \sin(\theta) (\sin^2(\theta) + 4 \cos^2(\theta))^{1/2} , 
\\
D & := (16 \sin^4(\theta) + 9 \cos^4(\theta) + 40 \sin^2(\theta) \cos^2(\theta))^{1/2} . 
\end{align*}
Thus, one has 
\begin{align*}
8 \, \Max K_{\sigma}^{(n>2)} = - 2 -3 \cos^2(\theta) + C , 
\quad 
8 \, \Max K_{\sigma}^{(n=2)} = - 2 - 3 \cos^2(\theta) + D . 
\end{align*}
Note that $C, D \ge 0$. 
By using 
\[
1 
= (\sin^2(\theta) + \cos^2(\theta))^2 
= \sin^4(\theta) + 2 \sin^2(\theta) \cos^2(\theta) + \cos^4(\theta) , 
\]
a straightforward calculation yields that 
\begin{align*}
C^2 - D^2 
& = 
6 \sin(\theta) 
((\sin^2(\theta) + 4 \cos^2(\theta))^{1/2} 
\\
& \qquad 
+ 9 
- 9 \cos^4(\theta) 
- 15 \sin^4(\theta) 
- 36 \sin^2(\theta) \cos^2(\theta))
\\
& = 
6 \sin(\theta) 
((1+3\cos^2(\theta))^{1/2} 
- \sin(\theta) (\sin^2(\theta) + 3 \cos^2(\theta))) 
\\
& = 
6 \sin(\theta) 
\left(
(1 + 3 \cos^2(\theta))^{1/2} 
- (1 + 3 \cos^2(\theta) - 4 \cos^6(\theta))^{1/2} 
\right) 
\\
& \ge 0 . 
\end{align*}
The equality holds if and only if $\sin(\theta) = 0$ or $\cos(\theta) = 0$. 
\end{proof}

As a result, for the homogeneous ruled minimal hypersurface $S(0)$, 
the maximum of the sectional curvature is $- 1/4$. 
This implies that $S(\theta)$ has negative sectional curvatures for small $\theta$.

\bibliographystyle{amsplain}

\end{document}